\documentclass{amsart}

\usepackage{amssymb}
\usepackage{amsfonts}
\usepackage{amsmath}
\usepackage{lastpage}
\usepackage[utf8]{inputenc}
\usepackage[legalpaper,bookmarks=true,colorlinks=true,linkcolor=blue,citecolor=blue]{hyperref}
\usepackage{graphicx}%
\usepackage{fancyhdr}
\usepackage{color}
\usepackage[mathlines]{lineno}
\usepackage{lscape}
\usepackage{epsfig}
\usepackage{natbib}
\usepackage{booktabs}
\usepackage{geometry}
\usepackage{tgbonum}
\fontfamily{qcr}\selectfont

\newtheorem{theorem}{Theorem}
\theoremstyle{plain}

\newtheorem{corollary}{Corollary}

\newtheorem{definition}{Definition}

\newtheorem{proposition}{Proposition}

\numberwithin{equation}{section}

\begin{document}
\Large
\title{A parallel treatment of semi-continuous functions with left and right-functions and some application in pedagogy}
\author{Gane Samb Lo}

\maketitle

\begin{abstract}  Left and right-continuous functions play an important role in Real analysis, especially in Measure Theory and Integration on the real line and in Stochastic processes indexed by a continuous real time. Semi-continuous functions are also of major interest in the same way. This paper aims at presenting a useful handling of semi-continuous function in parallel with the way left or right continuous functions are treated. For example, a lower or upper continuous function shares the property that it is measurable if it fails to be upper or lower continuous at most at a number of countable points with  a left or right-continuous when it fails to be left or right-continuous at most at a number of countable points. As a final result, the comparison between the Riemann and the Lebesgue integrals on compact sets is done in a very comfortable and comprehensive way. The frame used here is open to more further sophistication.\medskip

\noindent \textbf{R\'esum\'e}.  Les fonctions continues \`a droite et/ou \`a gauche jouent un r\^ole important en Analyse r\'eelle, surtout en Th\'eorie de la Mesure et de l'Int\'egration, et en processus stochastiques index\'es par un temps continu r\'eel. Les fonctions semi-continues aussi jouent un r\^ole tout aussi important. Ce papier est l'occasion d'offrir une pr\'esentation des fonctions semi-continues parall\`element aux propri\'et\'es des fonctions continues à  droite et/ou \`a  gauche. Par exemple, une fonction semi-continue inf\'erieurement ou sup\'erieurement partage la propri\'et\'e qu'elle est encore mesurable si elle est semi-continue inf\'erieurement ou sup\'erieurement en dehors d'un ensemble au plus d\'enombrable de points. Finalement, l'approche utilis\'ee rend plus compr\'ehensible et plus confortable la comparaison entre les int\'egrales de Riemman et de Lebesgue sur un intervalle compact de la droite r\'eelle. L'approche est ouverte \`a  de futurs 
d\'eveloppements.\medskip

\noindent Gane Samb Lo : gane-samb.lo@ugb.edu.sn, gslo@aust.edu.ng, ganesamblo@ganesamblo.net.\\
LERSTAD, Gaston Berger University, Saint-Louis, Sénégal.\\
Affiliated to LSTA, Université Pierre et Marie Curie, Paris, France\\
Associated with African University of Sciences and Technology, AUST, Abuja, Nigeria\\
1178 Evanston Drive NW, T3P 0J9, Calgary, Alberta, Canada.\\

\end{abstract}

\section{Introduction}
Left and right continuous functions play an important role in real analysis,
in particular in measure theory on the real line. In this paper, we proceed
to a study of semi-continuous functions following the way we do with let or
continuous functions. Though the material used here is classical, the
present treatment lead to a very comfortable handling of the mentioned
materials and their use in a number of problems of modern analysis.\newline

\noindent The ideas summarized here come from our teaching of measure theory
and integration and are closely related to Borel functions and the Lebesgue
and Riemann integration. One will feel the Baire's ideas behind all parts of
the text. We hope that this note might be considered as an input
in the classical measure theory teaching. This note might be written in a
compact and short form. But we rather give details to allow students to read
it.\newline

\noindent To define semi-continuous numerical functions, recall the superior
and inferior limit of a function $f:$ $\mathbb{R}\rightarrow \mathbb{R}$ at a
point $x\in \mathbb{R}$ : 
\begin{equation*}
f^{\ast }(x)=\limsup_{y\rightarrow x}f(y)=\lim_{\varepsilon \downarrow
0}\sup \{f(y),y\in ]x-\varepsilon ,x+\varepsilon \lbrack \}
\end{equation*}

\noindent and 
\begin{equation*}
f_{\ast }(x)=\liminf_{y\rightarrow x}f(y)=\lim_{\varepsilon \downarrow
0}\inf \{f(y),y\in ]x-\varepsilon ,x+\varepsilon \lbrack \}
\end{equation*}

\noindent As well we may define left superior and inferior of a function $f:$
$\mathbb{R}\rightarrow \mathbb{R}$ at a point $x\in \mathbb{R}$ : 
\begin{equation*}
f^{\ast ,\ell }(x)=\limsup_{y\rightarrow x}f(y)=\lim_{\varepsilon \downarrow
0}\sup \{f(y),y\in \rbrack x-\varepsilon ,x\lbrack \}
\end{equation*}

\noindent and 
\begin{equation*}
f_{\ast ,\ell }(x)=\liminf_{y\rightarrow x}f(y)=\lim_{\varepsilon \downarrow
0}\inf \{f(y),y\in \rbrack x-\varepsilon ,x\lbrack \}
\end{equation*}

\noindent and the right superior and inferior of a function $f:$ $\mathbb{R}%
\rightarrow \mathbb{R}$ at a point $x\in \mathbb{R}$ : 
\begin{equation*}
f^{\ast ,r}(x)=\limsup_{y\rightarrow x}f(y)=\lim_{\varepsilon \downarrow
0}\sup \{f(y),y\in \rbrack x ,x+\varepsilon \lbrack \}
\end{equation*}

\noindent and 
\begin{equation*}
f_{\ast ,r}(x)=\liminf_{y\rightarrow x}f(y)=\lim_{\varepsilon \downarrow
0}\inf \{f(y),y\in ]x,x+\varepsilon \lbrack \}
\end{equation*}

\bigskip \noindent These definitions are possible because the extrema we used in them
are either non-decreasing or non-increasing as $\varepsilon \downarrow 0$ and
the limits as $\varepsilon \downarrow 0$ are infima or suprema in \ $%
\overline{\mathbb{R}}$. For example, for the definition of $f^{\ast }(x)$
and $f_{\ast }(x),\sup \{f(y),y\in ]x-\varepsilon ,x+\varepsilon \}$ is
non-increasing as $\varepsilon \downarrow 0$ and $\inf \{f(y),y\in
]x-\varepsilon ,x+\varepsilon \}$ is non-decreasing $\varepsilon \downarrow
0.$\newline

\noindent The following properties also also known, whenever the
expressions make sense :\\

\noindent (1) $f$ is respectively continuous at $x\in \mathbb{R}$, or $f$ is
right-continuous, or $f $ is right continuous at $x$ if and only if, we
respectively have $\ f^{\ast }(x)=f_{\ast }(x),$ $\ $or $f^{\ast ,\ell
}(x)=f_{\ast ,\ell }(x),$ or $f^{\ast ,r}(x)=f_{\ast ,r}(x)$.\newline

\noindent (2) For any $x\in ,$ \ $(-f(x))^{\ast }=-f_{\ast }(x),$ \ $%
(-f(x))^{\ast ,\ell }=-f_{\ast ,\ell }(x),$\ $(-f(x))^{\ast ,r}=-f_{\ast
,r}(x).$\newline

\noindent In the sequel the functions $f^{\ast },f_{\ast },$ $f^{\ast ,\ell
},$ $f_{\ast ,\ell },$ $f^{\ast ,r}$ and $f_{\ast ,r}$ \ are defined on a
set $I\subset \mathbb{R}$, whenever they exist, by the following graphs : 
\begin{equation*}
x\hookrightarrow f^{\ast }(x),\text{ }x\hookrightarrow f_{\ast }(x),\text{ }%
x\hookrightarrow f^{\ast ,\ell }(x),x\hookrightarrow f_{\ast ,\ell
}(x),x\hookrightarrow f^{\ast ,r}(x).
\end{equation*}

\bigskip \noindent We have the following definitions.

\begin{definition}
We define\newline

\noindent \textbf{(1)} A function $f:\mathbb{R}\mathbb{\mapsto }\mathbb{R}$
is upper semi-continuous (usc) on an open interval $I$ of $\mathbb{R}$ if
and only if, for any $x\in I,$ $f(x)=f^{\ast }(x).$ It is lower
semi-continuous on $I$\ (lsc) if and ony if\thinspace , for any $x\in I,$ $%
f(x)=f_{\ast }(x).$\newline

\noindent \textbf{(2)} A function $f:\mathbb{R}\mathbb{\mapsto }\mathbb{R}$
is upper left-semi-continuous (ulsc) on an open interval $I$ of $\mathbb{R}$
\ if and only if, for any $x\in I,$ $f(x)=f^{\ast ,\ell }(x).$ It is lower
left-semi-continuous (llsc) on $I$ if and only if, for any $x\in I,$ $%
f(x)=f_{\ast ,\ell }(x).$\newline

\noindent (1) A function $f:\mathbb{R}\mathbb{\mapsto }\mathbb{R}$ is upper
right-semi-continuous (ursc) on an open interval $I$ of $\mathbb{R}$ if and
only if , for any $x\in I,$ $f(x)=f^{\ast ,r}(x).$ It is lower
right-semi-continuous (lrdc) on $I$ if and only if, for any $x\in I,$ $%
f(x)=f_{\ast ,r}(x).$\newline
\end{definition}

\noindent \textbf{Immediate remark}. A right-continuous function needs to
be lower right-semi-continuous and upper right-semi-continuous. But a lower
or upper semi-continuous function is not necessarily right or
left-continuous. So using \textit{ulsc}, \textit{ursc}, \textit{llsc} or\textit{\
lrdc} functions may extend the working classes in Analysis beyond
continuous or left or right-continuous functions. \newline

\noindent The definitions below are automatically extensible to numerical functions
defined on metric spaces. But let just see that we still have the classical
definitions of semi-continuous functions on $\mathbb{R}$. We have

\begin{proposition}
A function $f:\mathbb{R}\mathbb{\mapsto }\mathbb{R}$ is upper
semi-continuous (usc) at $x\in I,$ where $I$\ is an open interval of $%
\mathbb{R}$ if and only if

\begin{equation}
\forall \varepsilon >0,\text{ }\exists \eta >0,\text{ }(\left\vert
y-x\right\vert <\eta \text{ and }y\in I)\Longrightarrow f(y)\leq
f(x)+\varepsilon   \label{defclass}
\end{equation}%
A function $f:\mathbb{R}\mathbb{\mapsto }\mathbb{R}$ is lower
semi-continuous (usc) at $x\in I,$ where $I$\ is an open interval of $%
\mathbb{R}$ if and only if

\begin{equation*}
\forall \varepsilon >0,\text{ }\exists \eta >0,\text{ }(\left\vert
y-x\right\vert <\eta \text{ and }y\in I)\Longrightarrow f(x)+\varepsilon
\leq f(y).
\end{equation*}
\end{proposition}

\bigskip \noindent \textbf{Proof}. It is enough to prove only the first assertion,
since the second may be derived from the first and vice-versa. So, we have
to prove that (\ref{defclass}) is equivalent to $f^{\ast }(x)=f(x).$ Suppose
that (\ref{defclass}) holds. Then for any $\varepsilon >0,$ there exists $%
\eta >0$ such that $y\in ]x-\eta ,x+\eta \lbrack $ implies that $f(y)\leq
f(x)+\varepsilon .$ This means that for any $0<\delta \leq \eta ,$ 
\begin{equation*}
f(x)\leq f_{\delta }^{\ast }(x)=\sup \{]x-\delta ,x+\delta \lbrack \}\leq
f(x)+\varepsilon .
\end{equation*}

\bigskip \noindent We get, as $\delta \downarrow 0,$ for any $\varepsilon >0,$%
\begin{equation*}
f(x)\leq f^{\ast }(x)\leq f(x)+\varepsilon .
\end{equation*}

\bigskip \noindent We get $f^{\ast }(x)=f(x)$ by letting $\varepsilon \downarrow 0.$ Let us
prove the reverse implication, \textit{ad absordium}. Suppose that (\ref%
{defclass}) is false. Then there exists $\varepsilon >0$ such that for all
for all $\eta >0$, there exists $y(\eta )\in y\in ]x-\eta ,x+\eta \lbrack $
such that  $f(y(\eta ))>f(x)+\varepsilon ,$ which means that   
\begin{equation*}
f_{\eta }^{\ast }(x)>f(x)+\varepsilon .
\end{equation*}%
By letting $\eta \downarrow 0,$ we get%
\begin{equation*}
f^{\ast }(x)\geq f(x)+\varepsilon ,
\end{equation*}

\bigskip \noindent so that the equality $f^{\ast }(x)=f(x)$ does not holds.\\

\bigskip \noindent We are going to see that these functions have similar properties with right
and left-continuous functions.

\section{Semi-continuity and measurability}
\noindent In this section the most useful result is of the
following type

\begin{theorem}
\label{demiF_theo1} For any function $f$, $f^{\ast }$ is upper
semi-continuous and $f_{\ast }$ is lower semi-continuous.
\end{theorem}

\noindent From this result, we get the following proposition.

\begin{proposition}
\label{demiF_prop1} Upper and lower semi-continuous functions are measurable.
\end{proposition}

\bigskip \noindent \textbf{Proof of Proposition \ref{demiF_prop1}}. Let $f$
be upper semi-continuous. Set for each $n\geq 1,$%
\begin{equation*}
D_{n}=\{k2^{-n},k\in \mathbb{Z},n\geq 1\}.
\end{equation*}

\bigskip \noindent and

\begin{equation*}
f_{n}(x)=\sum_{k=-\infty }^{k=+\infty }\sup \{f(z),z\in
]_{]k2^{-n},(k+1)2^{-n}[}\}1_{]k2^{-n},(k+1)2^{-n}[}(x)+\sum_{s\in
D_{n}}f(s)1_{\{s\}}(x).
\end{equation*}

\bigskip \noindent For each $n\geq 1,f_{n}$ is measurable. Let us show that as $%
n\rightarrow +\infty ,$ for each $x\in \mathbb{R},$ $f_{n}(x)\rightarrow
f(x) $. We have two cases.\newline

\noindent If 
\begin{equation*}
x\in D=\bigcup_{n\geq 1} D_n,
\end{equation*}

\noindent we have, as $n\rightarrow +\infty$, 
\begin{equation*}
f_n(x)=f(x) \rightarrow f(x).
\end{equation*}

\noindent If

\begin{equation*}
x\notin D=\bigcup_{n\geq 1} D_n,
\end{equation*}

\bigskip \noindent then for $n\geq 1$, there exists $k=k(n,x)$ such that $%
x\in]k2^{-n},(k+1)2^{-n}[$. Then there exists $\varepsilon(n)>0$ such that
for any $0<\varepsilon\leq \varepsilon(n)$, we have

\begin{equation*}
]x-\varepsilon,x+\varepsilon[ \subset x\in]k(n,x)2^{-n},(k(n,x)+1)2^{-n}[, 
\text{ (A) }
\end{equation*}

\bigskip \noindent and next,

\begin{equation*}
\sup \{f(z),z\in ]x-\varepsilon,x+\varepsilon[\} \leq \sup \{f(z),z\in
[k2^{-n},(k+1)2^{-n}[\}. \text{ (B) }
\end{equation*}

\bigskip \noindent As well, for any $\eta>0$, for $2^{-n}<\eta$, for $%
\varepsilon(n)>0 $ such that (A) holds, for any $n\geq 1$, for any

\begin{equation*}
x\in]k(n,x)2^{-n},(k(n,x)+1)2^{-n}[ \in ]x-\eta,x+\eta[,
\end{equation*}

\noindent we have

\begin{equation*}
\sup \{f(z), z\in ]k(n,x)2^{-n},(k(n,x)+1)2^{-n}[\} \subset \sup \{f(z), z
\in ]x-\eta,x+\eta[ \}. \text{ (C) }
\end{equation*}

\bigskip \noindent By combining (A) and (B), we have for $\eta>0$, $2^{-n}<\eta$, for 
$\varepsilon(n)>0$ such that (A) holds, for any $0<\varepsilon\leq
\varepsilon(n)$, 

\begin{equation*}
f^{*,\varepsilon} \leq f_n \leq f^{*,\eta}.
\end{equation*}

\noindent Now let $\varepsilon\downarrow 0$ and only after, let $%
n\rightarrow +\infty$ to have

\begin{equation*}
f^{*} \leq \liminf_{n\rightarrow +\infty} f_n(x) \leq \limsup_{n\rightarrow
+\infty} f_n(x) \leq f^{*,\eta}.
\end{equation*}

\noindent Finally, let $\varepsilon\downarrow 0$ to get

\begin{equation*}
f_n(x)\rightarrow f^{*}=f(x).
\end{equation*}

\noindent Then, we have $f_{n}(x)=f(x)\rightarrow f(x)$ as $n\rightarrow
+\infty $. We conclude that $f$ is measurable.\newline

\noindent This proves the measurability of an upper semi-continuous. To
prove this for a lower semi-continuous, use the relation that for an lower
semi-continuous $f$, $-f$ is an upper semi-continuous and the conclusion is
immediate.$\square$\newline

\noindent Before we go further, let us mention that the same proof
is valid for the following corollary.

\begin{corollary}
\label{demiF_cor01} Let $f$ be real-valued function defined on $I=[a,b]$ or $%
I=\mathbb{R}$ and for each $n$, let $I_{i,n}=(x_{i,n},x_{i+1,n})$, $i\in
J(n) $, be consecutive intervals, with non zero length, which partitions $I$
such that

\begin{equation*}
\sup_{i\in J(n)}|x_{i+1,n}-x_{i+1,n}|\rightarrow 0\ as\ n\rightarrow +\infty.
\end{equation*}

\noindent Denote $D=\cup _{n}\cup _{i\in J(n)}\{x_{i,n}\}$. We have for any $%
x\notin D$, 
\begin{equation*}
H_{n}(x)=\sum_{i}\sup \{f(z),z\in
]x_{i,n},x_{i+1,n}[\}1_{]x_{i,n},x_{i+1,n}[}(x)\rightarrow f^{\ast }(x).
\end{equation*}

\bigskip \noindent and

\begin{equation*}
h_{n}(x)=\sum_{i}\inf \{f(z),z\in
]x_{i,n},x_{i+1,n}[\}1_{]x_{i,n},x_{i+1,n}[}(x)\rightarrow f_{\ast }(x).
\end{equation*}
\end{corollary}

\bigskip \noindent \textbf{Proof of Theorem \ref{demiF_theo1}}. Set $%
g=f^{\ast }$. Let us show that $g^{\ast }=g.$ Let $x$ be fixed. By
definition for $\varepsilon >0, $

\begin{equation*}
g_{\varepsilon }^{\ast }(x)=\sup \{g(z),z\in ]x-\varepsilon ,x+\varepsilon
\lbrack \}=\sup \{f^{\ast }(z),z\in ]x-\varepsilon ,x+\varepsilon \lbrack \}
\end{equation*}

\bigskip \noindent Suppose that $g^{\ast }(x)$ is finite. So is $g_{\varepsilon
}^{\ast }(x)$ for $\varepsilon $ small enough. Use the characterization of
the supremum : for all $\eta >0$, there exists $z_{0}\in ]x-\varepsilon
,x+\varepsilon \lbrack $ such that

\begin{equation*}
g_{\varepsilon }^{\ast }(x)-\eta <f^{\ast }(z_{0})<g_{\varepsilon }^{\ast
}(x).
\end{equation*}

\noindent Since $z_{0}\in ]x-\varepsilon ,x+\varepsilon \lbrack ,$ there
exists $\varepsilon _{0}>0$ such that $z_{0}\in ]z_{0}-\varepsilon
_{0},z_{0}+\varepsilon _{0}[\subset z]x-\varepsilon ,x+\varepsilon \lbrack .$
Remind that $f_{h}^{\ast }(z_{0})\downarrow f^{\ast }(z_{0})$ as $%
h\downarrow 0.$ Then for $0<h\leq \varepsilon _{0},$

\begin{equation*}
f_{h}^{\ast }(z_{0})=\sup \{f(z),z\in ]z_{0}-h,z_{0}+h[\}\leq \sup
\{f(z),z\in ]x-\varepsilon ,x+\varepsilon \lbrack \}=f_{\varepsilon }^{\ast
}(x)
\end{equation*}

\bigskip \noindent so that 
\begin{equation*}
g_{\varepsilon }^{\ast }(x)-\eta <f^{\ast }(x).
\end{equation*}

\noindent Let $\eta \downarrow 0$ and next $\varepsilon \downarrow 0$ to get%
\begin{equation}
g^{\ast }(x)\leq f^{\ast }(x).  \tag{INEQ}
\end{equation}

\noindent Suppose that these inequality is strict, that is 
\begin{equation*}
g^{\ast }(x)<f^{\ast }(x).
\end{equation*}

\noindent So, we can find $\eta >0$ such that

\begin{equation*}
g^{\ast }(x)<f^{\ast }(x)-\eta .
\end{equation*}

\noindent Since $g_{\varepsilon }^{\ast }(x)\downarrow g^{\ast }(x)$, there
exists $\varepsilon _{0}$ such that 
\begin{equation*}
g^{\ast }(x)\leq g_{\varepsilon _{0}}^{\ast }(x)<f^{\ast }(x)-\eta,
\end{equation*}

\noindent that is 
\begin{equation*}
g^{\ast }(x)\leq \sup \{g(z),z\in ]x-\varepsilon _{0},x+\varepsilon
_{0}[\}<f^{\ast }(x)-\eta,
\end{equation*}

\noindent that is also

\begin{equation*}
g^{\ast }(x)\leq \sup \{f^{\ast }(z),z\in ]x-\varepsilon _{0},x+\varepsilon
_{0}[\}<f^{\ast }(x)-\eta .
\end{equation*}

\noindent But 
\begin{equation*}
x\in ]x-\varepsilon _{0},]x+\varepsilon _{0}[
\end{equation*}

\noindent and then 
\begin{equation*}
f^{\ast }(x)\leq \sup \{f^{\ast }(z),z\in ]x-\varepsilon _{0},]x+\varepsilon
_{0}[\}<f^{\ast }(x)-\eta .
\end{equation*}

\bigskip \noindent We arrive at the conclusion 
\begin{equation*}
f^{\ast }(x)<f^{\ast }(x)-\eta .
\end{equation*}

\bigskip \noindent This is absurd. Hence the inequality (INEQ) is an equality, that is%
\begin{equation*}
g^{\ast }(x)=f^{\ast }(x).
\end{equation*}

\noindent Now suppose that $g^{\ast }(x)$ is infinite. This means that for
any $A>0,$ we can find in any interval $]x-\varepsilon ,x-\varepsilon\lbrack$, a point $z(\varepsilon )$ such that 
\begin{equation*}
g(z(\varepsilon ))>A,
\end{equation*}

\noindent Now fix $A$ and $\varepsilon >0$. Consider $z(\varepsilon )$ such
that the last inequality holds. Since $z(\varepsilon )\in ]x-\varepsilon
,]x-\varepsilon \lbrack ,$ there exists $r_{0}>0$ such that%
\begin{equation*}
z(\varepsilon )\in ]z(\varepsilon )-r_{0},z(\varepsilon )+r_{0}[\subset
]x-\varepsilon ,x+\varepsilon \lbrack .
\end{equation*}

\noindent So, since 
\begin{equation*}
f_{r}^{\ast }(z(\varepsilon ))\downarrow g(z(\varepsilon ))=f_{r}^{\ast
}(z(\varepsilon ))>A,
\end{equation*}

\noindent then for $r$ small enough%
\begin{equation*}
f_{r}^{\ast }(z(\varepsilon ))>A/2,
\end{equation*}

\noindent that is

\begin{equation*}
\sup \{f(z),z\in ]z(\varepsilon )-r,z(\varepsilon )+r[\}>A/2
\end{equation*}

\noindent and again there will be a $u\in z\in ]z(\varepsilon
)-r,z(\varepsilon )+r[$ such that

\begin{equation*}
f(u)>A/4.
\end{equation*}

\noindent If $r$ is taken small enough such that $r\leq r_{0}$, we have
found a point $u\in ]x-\varepsilon ,]x-\varepsilon \lbrack $ such that%
\begin{equation*}
f(u)>A/4.
\end{equation*}

\noindent We have proved that for any $A>0,$ for any $\varepsilon >0,$ we
can find a point in $u\in ]x-\varepsilon ,]x-\varepsilon \lbrack $ such that

\begin{equation*}
f(u)>A/4
\end{equation*}

\noindent and hence

\begin{equation*}
f_{\varepsilon }^{\ast }(x)=\sup \{f^{\ast }(z),z\in ]x-\varepsilon
,]x+\varepsilon \lbrack \}\geq A/4.
\end{equation*}

\bigskip \noindent Let $\varepsilon \downarrow 0$ and next $A\uparrow
+\infty $ to get that%
\begin{equation*}
g(x)=f^{\ast }(x)=+\infty .
\end{equation*}

\noindent Thus 
\begin{equation*}
g^{\ast }(x)=g(x)=+\infty .
\end{equation*}

\bigskip \noindent We have finished and we proved that%
\begin{equation*}
g^{\ast }=g.
\end{equation*}

\noindent \textbf{Conclusion} : if $g^{\ast }$ is finite, then it is upper
semi-continuous. We similarly use that if $g_{\ast }$ is finite, then it is
lower semi-continuous , by exploiting the formula $(-f)^{\ast }=-f_{\ast }.$%
\newline

\bigskip \noindent We get the following useful result

\begin{corollary}
\label{demiF_cor02} \label{cor01} Any function $f$ which is lower
semi-continuous except at a countable number of points and any function $g$
which is upper semi-continuous except at a countable number points are
measurable.
\end{corollary}

\bigskip \noindent \textbf{Proof}. We give the proof only for a lower
semi-continuous function $f$. The other case is derived from the application
of the results of the first case to the opposite of $f$.\newline

\noindent Let us denote the countable $D$ set on which $f$ is not lower
semi-continuous. We have 
\begin{equation*}
f=f_{\ast }\text{ on }D^{c}
\end{equation*}

\noindent and then 
\begin{equation}
f=f_{\ast }1_{D^{c}}+\sum_{s\in D}f(s)1_{\{s\}}.  \label{decomFlscN}
\end{equation}

\noindent By Theorem \ref{demiF_theo1} and \ref{demiF_prop1}, $f_{\ast }$ is
lower semi-continuous. Hence by (\ref{decomFlscN}), $f$ is measurable.%
\newline

\subsection{Application to a rephrasing of the classical comparison between
Riemann and Lebesgue integrals}

\begin{proposition}
\label{demiF_prop02} Let $[a,b]$ be a bounded and nonempty compact interval
of $\mathbb{R}$ and $\lambda$ be the Lebesgue measure on $[a,b]$. Let $%
f:[a,b] \mapsto \mathbb{R}$ be a bounded function by $M$. Then $f$ is
Riemann integrable if and only if $f$ is almost $\lambda$-a.e. continuous,
and its Riemann integral is equal to its Lebesgue integral.
\end{proposition}

\noindent \textbf{Proof of Proposition \ref{demiF_prop02}}. In that proof, for all rules using measure theory and integration, we refer to \cite{ips-mestuto-ang}. Consider the finite non-empty and bounded interval $]a,b]$ of $%
\mathbb{R}$ and a bounded, say by $M$, real-valued function $f$ defined on $]a,b]$. For each $n\geq 1$, consider a subdivision of $]a,b]$, which
partitions $]a,b]$ into the $\ell (n)$ sub-intervals, in the form,

\begin{equation*}
]a,b]=\sum_{i=0}^{\ell (n)-1}]x_{i,n},x_{i+1,n}].
\end{equation*}

\bigskip \noindent with modulus 
\begin{equation*}
m(\pi _{n})=\max_{0\leq i\leq \ell(n)}(x_{i+1,n}-x_{i+1,n}) \rightarrow 0 \ as
\ n\rightarrow +\infty.
\end{equation*}

\bigskip \noindent Define for each $n\geq 0$, for each $i$, $0\leq i \leq \ell(n)-1$

\begin{equation*}
m_{i,n}=\inf \{f(z),x_{i,n}\leq z<x_{i+1,n}\}\text{ and }M_{i,n}=\sup
\{f(z),x_{i,n}\leq z<x_{i+1,n}\}
\end{equation*}

\begin{equation*}
h_{n}=\sum_{i=0}^{\ell (n)-1}m_{i,n}1_{]x_{i,n},x_{i+1,n}]}\text{ and }%
H_{n}=\sum_{i=0}^{\ell (n)-1}M_{i,n}1_{]x_{i,n},x_{i+1,n}]}.
\end{equation*}

\noindent and

\begin{equation*}
D=\bigcup_{k}\{x_{0,n},x_{1,n},...,x_{\ell (n),n}\}.
\end{equation*}

\noindent Remark that $f^{\ast }$ and $f_{\ast }$ are bounded by $M$ and $D$
is countable. For a fixed $x$ in $x\in \lbrack a,b]\setminus D$, for any $%
n\geq 1$, there exists $i$, $0\leq i(n)\leq \ell (n)-1$\ such that $x\in
]x_{i(n),n},x_{i(n)+1,n}[.$ Then, we can find $\varepsilon >0$\ such that 
\begin{equation*}
]x-\varepsilon ,x+\varepsilon \lbrack \subseteq ]x_{i(n),n},x_{i(n)+1,n}[.
\end{equation*}

\noindent It is clear that $x\in \lbrack a,b]\setminus D$

\begin{equation*}
h_{n}(x)\leq f_{\ast }^{\varepsilon }(x).
\end{equation*}

\noindent and

\begin{equation*}
H_{n}(x) \geq f_{\ast }^{\varepsilon }(x).
\end{equation*}

\noindent By letting $\varepsilon \downarrow 0$, we get

\begin{equation}
h_{n}(x)\leq f_{\ast }(x)\leq f^{\ast }\leq H_{n}(x).  \label{B1}
\end{equation}

\bigskip \noindent So (\ref{B1}) holds for any $x\in \lbrack a,b]\setminus D$. Now,
let $\eta >0$. For any fixed $n\geq 1$, we may use the characterization of the
suprema and the infima on $\mathbb{R}$, and show that the Lebesgue integral
of $h_{n}$ and $H_{n} $ can be approximated to Riemann sums in the form

\begin{equation*}
\int H_{n}d\lambda =\sum_{i=0}^{\ell
(n)-1}M_{i,n}(x_{i+1,n}-x_{i,n})<\sum_{i=0}^{\ell
(n)-1}f(d_{i,n})(x_{i+1,n}-x_{i,n})+\eta
\end{equation*}

\noindent and

\begin{equation*}
\int h_{n}d\lambda =\sum_{i=0}^{\ell
(n)-1}m_{i,n}(x_{i+1,n}-x_{i,n})>\sum_{i=0}^{\ell
(n)-1}f(c_{i,n})(x_{i+1,n}-x_{i,n})-\eta.
\end{equation*}

\noindent By putting together the previous facts, we get

\begin{eqnarray*}
&&\sum_{i=0}^{\ell (n)-1}f(c_{i,n})\left(x_{i+1,n}-x_{i,n}\right)
-\eta \\
&\leq &\int_{[a,b]\setminus D}f_{\ast }(x)d\lambda \\
&\leq &\int_{[a,b]\setminus D}f^{\ast }(x)d\lambda \\
&\leq &\sum_{i=0}^{\ell (n)-1}f(d_{i,n})(x_{i+1,n}-x_{i,n})+\eta,
\end{eqnarray*}

\noindent where we use that $D$ is countable and then, is a $\lambda $-null
set. By letting $\eta \rightarrow 0$, we get

\begin{equation*}
\int_{\lbrack a,b]\setminus D}f_{\ast}(x) \ d\lambda(x)=\int_{[a,b]\setminus D}f^{\ast
}(x) \ d\lambda(x).
\end{equation*}

\bigskip \noindent Thus, we have that

\begin{equation*}
f_{\ast }=f^{\ast }=f
\end{equation*}

\bigskip \noindent outside a countable subset $[a,b]\setminus D,$ then outside a $%
\lambda$-null subset of $[a,b]$ and 
\begin{equation}
\int_{a}^{b}f(x)\ d\lambda(x)=\int_{[a,b]}f(x)d\lambda (x).  \label{riemann-lebesgues}
\end{equation}

\noindent Now suppose that $f$ is $a.e.$ continuous, that is $f_{\ast
}=f^{\ast }=f$ $\lambda -a.e.$ Consider a sequence of Riemann sums for which
the sequence of modulus tends to zero with $n$ :

\begin{equation*}
S_{n}=\sum_{i=0}^{\ell (n)-1}f(c_{i,n})(F(x_{i+1,n})-F(x_{i,n})).
\end{equation*}

\noindent Let us denote the $\lambda$-null set $H=(f^{\star}\neq f_{\star})$We have

\begin{equation*}
\int h_{n}d\lambda \leq S_n \leq \int H_{n}d\lambda
\end{equation*}

\noindent Denoting $H=[a,b]\setminus D,$ we have

\begin{equation*}
\int h_{n} 1_{H} d\lambda \leq S_n \leq \int H_{n} h_{n} 1_{H} d\lambda
\end{equation*}

\noindent By Corollary \ref{demiF_cor01} above, $h_{n}1_{H}\rightarrow
f_{\ast }1_{H}$ and $H_{n}1_{H}\rightarrow f^{\ast}1_{H}$ as $n\rightarrow
+\infty$. Since $|f|\leq M$ which is $\lambda$-integrable on $[a,,b]$, we may apply the Dominated Convergence Theorem to get

\begin{equation*}
\int f_{\ast }1_{H}d\lambda \leq \liminf_{n\rightarrow +\infty }S_{n}\leq
\limsup_{n\rightarrow +\infty }S_{n}\leq \int f_{\ast }1_{H}d\lambda .
\end{equation*}

\noindent From there, we have

\begin{equation*}
\liminf_{n\rightarrow +\infty }S_{n}=\limsup_{n\rightarrow +\infty
}S_{n}=\int f_{\ast }1_{H}d\lambda =\int f_{\ast }1_{H}d\lambda =\int
fd\lambda .
\end{equation*}

\noindent So, all Riemann sums converge to 
\begin{equation*}
I=\int fd\lambda \in \mathbb{R},
\end{equation*}

\noindent when the modulus goes to zero. Then $f$ is Riemann integrable and Formula \ref{riemann-lebesgues} holds again.\\

\noindent \textbf{Acknowledgment}. The author wishes to express his warm thanks to AUST and to professor Charles E. Chidume for hospitality while finalizing this work. He also acknowledges continuous support of the World Bank CEA-MITIC and Gaston Berger University of Saint-Louis.

\end{document}